\documentclass[reqno]{amsart}
\usepackage{latexsym,amssymb,amsfonts}
\usepackage{graphicx,color,graphicx}
\usepackage[mathscr]{euscript}
\usepackage{lscape}
\newtheorem{lemma}{\bf Lemma}[section]
\newtheorem{proposition}{\bf Proposition}[section]
\newtheorem{corollary}{\bf Corollary}

\newtheorem{theorem}{\bf Theorem}[section]

\numberwithin{equation}{section} \theoremstyle{plain}
\theoremstyle{definition}

\input cyracc.def

\begin{document}
\title[On tail distributions of supremum and quadratic variation]
{On tail distributions of supremum and quadratic variation of
local martingales}

\author{Liptser R.}
\address{Dept. Electrical Engineering-Systems,
 Tel Aviv University, 69978 Tel Aviv, Israel}
\email{<liptser@eng.tau.ac.il>}
\author{Novikov A.}
\address{Dept. Mathematical Sciences, University of Technology Sydney, PO Box, 123.
Broadway, NSW 2007, Australia} \email{<prob@maths.uts.edu.au> }

\begin{abstract}
We extend some known results relating the distribution tails of a
continuous local martingale supremum and its quadratic variation
to the case of locally square integrable martingales with bounded
jumps. The predictable and optional quadratic variations are
involved in the main result.
\end{abstract}

\maketitle

\section{Introduction and main result}

Denote by $\mathscr{M}(\mathscr{M}_{\rm loc})$ and
$\mathscr{M}^2 (\mathscr{M}^2_{\rm loc} , \mathscr{M}^c_{\rm
loc})$ the classes of all martingales (local martingale) and
 square integrable
(locally square integrable, continuous local martingales)
$M=(M_t)_{t\ge 0},\, M_0=0$ (with paths in the Skorokhod space
$\mathbb{D}_{[0,\infty)}$) defined on
$(\Omega,\mathscr{F},(\mathscr{F}_t)_{t\ge 0}, P)$ a stochastic
basis with standard general conditions. Recall that any random
process $X$ with paths in the Skorokhod space and defined on the
above-mentioned stochastic basis belongs to the class
$\mathcal{D}$ if the family $(X_\tau,\tau\in\mathcal{T})$, where
$\mathcal{T}$ is the set of stopping times $\tau$, is uniformly
integrable.

Henceforth $\triangle M_t:=M_t-M_{t-}$, $\langle M\rangle_t$ and
$[M,M]_t$ are the jumps, predictable quadratic variation and optional
quadratic variation processes of $M$ respectively.

It is well-known (see e.g. \cite{LSMar}, \cite{JSn} and references therein)
that for local martingales from
$\mathscr{M}^2_{\rm loc}$:
$$
\langle M\rangle_\infty<\infty, \ {\rm a.s.}\Rightarrow
  \begin{cases}
[M,M]_\infty<\infty \  {\rm a.s.}
\\
\lim\limits_{t\to\infty}M_t=M_\infty\in \mathbb{R} \ {\rm a.s.}
  \end{cases}
$$
There are many other
 well-known relations between $M_\infty$ and $\langle M\rangle_\infty$
 (e.g., Burkholder--Gundy--Davis's
inequalities, law of large numbers for martingales, etc.) which
are valid for local martingales with jumps.

If $M\in \mathscr{M}\cap \mathcal{D}$, then
$M$ satisfies the Wald equality:
$$EM_\infty =0$$
which plays a fundamental role in many applications in stochastic
analysis. Often, a direct verification of the uniform
integrability is difficult. In this connection, we mention one
result from Novikov, \cite{Nov}, establishing a relation between
the tail distributions of $\langle M\rangle_\infty$ and
$EM_\infty$. A similar result is also proved in Elworthy, Li and
Yor, \cite{ELY}, under slightly different conditions than in
\cite{Nov}. Concerning the related topic dealing with a one-sided
stochastic boundary, see Pe\u{s}kir and Shiryaev, \cite{PeSh}, and
Vondra\u{c}ek \cite{Vod}.

\bigskip
\noindent {\bf Theorem*.} {\it Let $M\in\mathscr{M}^c_{\rm loc}$
\emph{and }$\langle M\rangle_\infty<\infty$ a.s. If
$\sup_{t>0}Ee^{\varepsilon M_t}<\infty$ for some positive
$\varepsilon$,
 then}\footnote{$a^+=\max(a,0),
a^-=\max(-a,0)$} $0\le EM_\infty \le EM^{+}_\infty<\infty$
\emph{and}
$$
\lim_{\lambda\to\infty}\lambda P\big(\langle
M\rangle^{1/2}_\infty> \lambda\big)=\sqrt{\frac{2}{\pi}}
EM_\infty.
$$
One of our goals is a generalization of Theorem* statement  for
local martingales with bounded jumps.

\begin{theorem}\label{theo-n1} Let $M\in\mathscr{M}^2_{\rm loc},
\ \langle
M\rangle_\infty<\infty$ a.s. and $M^+\in \mathcal{D}$. Then:\\
{\rm (i)} $\lim_{t\to\infty}M_t:=M_\infty$ exists and
$$0\le EM_\infty\le EM^+_\infty<\infty;$$
{\rm (ii)} $|\triangle M|\in \mathcal{D}$ and {\rm(}i{\rm)}
provide
$$
\lim\limits_{\lambda\to\infty}\lambda P\big(\sup_{t\ge 0} M^-_t>
\lambda\big)=EM_\infty;
$$
{\rm (iii)} $ |\triangle M|\le K $ and
\begin{equation}\label{oxbla0}
Ee^{\varepsilon M_\infty}<\infty,
\end{equation}
for some positive $K$ and $\varepsilon$, provide
$$ \lim\limits_{\lambda\to\infty}\lambda P\big(\langle
M\rangle^{1/2}_\infty>\lambda\big)
=\lim\limits_{\lambda\to\infty}\lambda P\big([M,M]^{1/2}_\infty>
\lambda\big)=\sqrt{\frac{2}{\pi}}EM_\infty.
$$
\end{theorem}

If $M^+\in \mathcal{D}$, Theorem \ref{theo-n1} gives necessary and
sufficient conditions for $M\in \mathcal{D}$ expressed in terms of
$\sup_{t\ge 0}M^-_t$, $\langle M\rangle_\infty$, and
$[M,M]_\infty$ which are useful in some applications (see, e.g.,
by  Jacod and Shiryaev \cite{JS98}).

\bigskip
\begin{corollary}\label{cor-0}
Under the assumptions of Theorem \ref{theo-n1}, the process $ M\in
\mathcal{D}$ iff any of the following conditions holds:
$$\lim\limits_{\lambda\to\infty}\lambda
P\big(\sup_{t\ge 0} M^-_t> \lambda\big)=0,
$$
$$
\lim\limits_{\lambda\to\infty}\lambda P\big(\langle
M\rangle^{1/2}_\infty>\lambda\big)=0,
$$
$$
\lim\limits_{\lambda\to\infty}\lambda P\big([M,M]^{1/2}_\infty>
\lambda\big)=0.
$$
\end{corollary}

A few publications preceded \cite{Nov} and \cite{ELY} (see Azema,
Gundy and Yor, \cite{AGY}; for discrete time martingales, Gundy,
\cite{ Gundy}, and Galtchouk and Novikov, \cite{GalNov}). Takaoka,
\cite{Takaoka}, presented a result similar to Theorem *.

The proofs of parts (i) and (ii) of Theorem \ref{theo-n1} are obvious
and  might even be known. The proof of (iii) exploits a
combination of techniques:
\begin{center}
``Stochastic exponential + Tauberian theorem''
\end{center}
which seems to have been firstly used by Novikov, \cite{Nov2}, to
obtain asymptotics  of the first passage times for Brownian
motion (see also \cite{Nov}) and for random walks (see, Novikov \cite{Nov3}).
Some necessary facts on the stochastic exponential are gathered in
Section \ref{section 2}. The proofs are given in Section \ref{section 3}.

The uniform boundedness assumption for $\Delta M$ might be
weakened by applying a standard "truncation" technique under some
additional assumptions on the tails distribution of $\Delta M$. We
show in Theorem \ref{theo-n3} that the uniform boundedness
assumption for $\Delta M$ is avoided if the stochastic
exponential possesses an evaluation in terms of $\langle
M\rangle_\infty$. This condition is borrowed from \cite{Nov} where
it is effectively applied for discrete-time martingales involving in a popular
gambling strategies.

\section{Preliminaries}
\label{section 2}
\subsection{Stochastic exponential}

For discontinuous martingales, the stochastic exponential has an
``intricate'' structure. So, we start with recalling the necessary
notions and objects involving in (ii)
(for more details, see e.g. \cite{LSMar} or \cite{JSn}).

For $M\in\mathscr{M}^2_{\rm loc}$, $M_0=0$, the decomposition
$M=M^c+M^d$ is well known, where $M^c,M^d\in\mathscr{M}^2_{\rm
loc}$ and are continuous and purely discontinuous
martingales respectively. Moreover, $\langle M\rangle= \langle
M^c\rangle+\langle M^d\rangle$, so the assumption $\langle
M\rangle_\infty<\infty$ provides $\langle
M^c\rangle_\infty<\infty$, $\langle M^d\rangle_\infty<\infty$. The
measure $\mu$ is associated with the jump process $\triangle
M\equiv \triangle M^d$ in the sense that for any measurable set
$A$ and $t>0$ $\mu((0,t]\times A)=\sum\limits_{s\le t}I(\triangle
M_s\in A)$. Denote by $\nu=\nu(dt,dz)$ its compensator. The
condition $|\triangle M|\le K$ provides the existence of a version
$\nu$ such that $\nu(\mathbb{R}_+\times\{|z|>K\})=0$. This version
of $\nu$ is used in the sequel.

The purely discontinuous martingale $M^d$ is defined as the It\^o
integral with respect to $\mu-\nu$:
$$
M^d_t=\int_0^t\int_{|z|\le K}z\big(\mu(ds,dz)-\nu(ds,dz)\big).
$$
Recall also that $ \int_{|z|\le K}z\nu(\{t\},dz)=0 \ {\rm a.s.} $
and
$$
\langle M^d\rangle_t= \int_0^t\int_{|z|\le K}z^2\nu(ds,dz)<\infty
\ {\rm a.s.}, \ t>0.
$$

Hence, $\langle M^d\rangle_t<\infty$ a.s. provides
\begin{equation}\label{2.1z}
\int_0^\infty\int_{|z|\le K}z^2\nu(ds,dz)<\infty \ {\rm a.s.}
\end{equation}

 This fact is important for further
considerations as long as we will deal with the cumulant process
$$
G_t(\lambda)=\int_0^t\int_{|z|\le K} \big(e^{\lambda z}-1-\lambda
z\big)\nu(ds,dz), \ \lambda\in\mathbb{R}.
$$
The boundedness of jumps and \eqref{2.1z} implies the existence of
$G_t(\lambda)$ and
$G_\infty(\lambda):=\lim_{t\to\infty}G_t(\lambda)<\infty$. The
cumulant process $G(\lambda)$, being increasing, possesses a
nonnegative jumps process
$$
\triangle G_t(\lambda):=\int_{|z|\le K}\big(e^{\lambda
z}-1-\lambda z\big) \nu(\{t\},dz).
$$

\medskip
A random process $\mathscr{E}(\lambda)$ with
\begin{equation}\label{2.2}
\mathscr{E}_t(\lambda)=\exp\Big(\frac{\lambda^2}{2}\langle
M^c\rangle_t+G_t(\lambda)\Big) \prod_{0<s\le t}\big(1+\triangle
G_s(\lambda)\big) e^{-\triangle G_s(\lambda)}
\end{equation}
is known as ``stochastic exponential'' for the martingale $M$. Note that
$\mathscr{E}_t>0$, since $\triangle G(\lambda)\ge 0$.

A remarkable property of the stochastic exponential is that the
process $\mathfrak{z}(\lambda)$,
\begin{equation}\label{A.2}
\mathfrak{z}_t(\lambda)=e^{\lambda M_t-\log\mathscr{E}_t(\lambda)}
\end{equation}
is a positive local martingale. Indeed, applying the It\^o formula
to \eqref{A.2}, we get
$$
d\mathfrak{z}_t(\lambda)=\lambda \mathfrak{z}_t(\lambda)
dM^c_t
+\int_{|z|\le K}\mathfrak{z}_{t-}(\lambda)
\frac{\big(e^{\lambda z}-1\big)}{1+\triangle
G_t(\lambda)}(\mu-\nu)(dt,dz),
$$
where the right-hand side is a sum of two local martingales. As any
nonnegative local martingale, $\mathfrak{z}(\lambda)$ is also a
supermartingale too (see e.g. Problem 1.4.4 in Liptser and Shiryaev \cite{LSMar}).
The
latter provides the existence of
$$
\mathfrak{z}_\infty(\lambda):=\lim_{t\to\infty}\mathfrak{z}_t(\lambda)\in\mathbb{R}_+
\ {\rm a.s.}
$$
with $E\mathfrak{z}_\tau(\lambda)\le 1$ for any Markov time $\tau$; hence,
in particular, $E\mathfrak{z}_\infty\le 1$.

\begin{proposition}\label{pro-2.1}
Let $|\triangle M|\le K$, $\langle M\rangle_\infty$ a.s. and condition
\eqref{oxbla0} hold. Then,
with $\varepsilon$ from \eqref{oxbla0} and any $\lambda\in
(0,\varepsilon]$,

{\rm 1)} $ E\mathfrak{z}_\infty(\lambda)=1$.

{\rm 2)}
$\mathscr{E}_\infty(\lambda)=\lim\limits_{t\to\infty}\mathscr{E}_t(\lambda)\in
\mathbb{R}_+$ a.s. and $\mathscr{E}_\infty(\lambda)>0$ a.s.
\end{proposition}
\begin{proof}
1) Let $(\tau_n)$ be an increasing sequence of stopping times,
$\lim_n\tau_n=\infty$, such that $(M_{t\wedge\tau_n})_{t\ge 0}$
and $(\mathfrak{z}_{t\wedge\tau_n}(\lambda))_{t\ge 0} \in
\mathcal{D}$ \ for any $n$. Then

\begin{equation}\label{=1}
E\mathfrak{z}_{\tau_n}(\lambda)\equiv 1.
\end{equation}
In order to finish the proof, we show that $\mathfrak{z}_{\tau_n}(\lambda)$
is majorized by uniformly integrable martingale
$
E\Big(e^{\lambda M^+_\infty}|\mathscr{F}_{\tau_{n}}\Big),
$
what is provided by \eqref{oxbla0}, applying Jensen's inequality:
$
E\big(e^{\lambda M^+_\infty}|\mathscr{F}_{\tau_{n}}\big) \ge
e^{\lambda  E(M^+_\infty|\mathscr{F}_{\tau_n})} \ge e^{\lambda
M^+_{\tau_n}}\ge \mathfrak{z}_{\tau_n}(\lambda).
$

Hence,
$(\mathfrak{z}_{\tau_n}(\lambda))_{n\ge 1} \in D$.

\medskip
2) Since
$
\mathfrak{z}_\infty(\lambda)=e^{\lambda M_\infty-\log\mathscr{E}_\infty(\lambda)}
$
with $\log 0=-\infty$, the desired  property holds true
provided that $\mathfrak{z}_\infty(\lambda)<\infty$ a.s.
\end{proof}

\section{The proof of Theorem \ref{theo-n1}}
\label{section 3}
\subsection{The proof of parts (i) and (ii)}

\mbox{}
 1) Let $(\tau_n)_{n\ge 1}$ be an increasing sequence of stopping times,
$\lim_n\tau_n=\infty$, such that $(M_{\tau_n})_ {n\ge 1} \in \mathcal{D}$
and, therefore,
$EM^-_{\tau_n}=EM^+_{\tau_n},n\ge 1$. Due to the assumption $M^+
\in \mathcal{D}$, we have
$\lim\limits_{n\to\infty}EM^+_{\tau_n}=EM^+_\infty<\infty$. Now,
applying the Fatou theorem, we find that $ EM^+_\infty\ge
EM^-_\infty$.

Hence,
$$
EM^+_\infty\ge EM^+_\infty-EM^-_\infty=EM_\infty\ge 0.
$$

\medskip
2) Set $S_\lambda=\inf\{t:M^-_{t}\ge\lambda\}$ and notice that
$$
\{S_\lambda<\infty\}=\{\sup_{t\ge 0} M^-_t>\lambda\}.
$$

Since $\triangle M_\infty=0$ and $|\triangle M| \in \mathcal{D}$,
the process  $(M_{t\wedge S_\lambda})_{t\ge 0}$ is a uniformly
integrable martingale with $EM_{S_\lambda}=0$.

Write
$$
\begin{aligned}
0=EM_{S_\lambda}&=EM_\infty
I_{\{S_\lambda=\infty\}}+EM_{S_\lambda}I_{\{S_\lambda<\infty\}}
\\
&=EM_\infty I_{\{S_\lambda=\infty\}}+E M_{S_\lambda}
I_{\{\sup_{t\ge 0}(-M_t)\ge\lambda\}}
\\
&=EM_\infty I_{\{S_\lambda=\infty\}}+E(M_{S_\lambda}-\lambda)
I_{\{S_\lambda<\infty\}}
\\
&\quad +\lambda P(\sup_{t\ge 0}M^-_t>\lambda).
\end{aligned}
$$

Finally, $EM^+_\infty<\infty$ provides
$\lim\limits_{\lambda\to\infty}S_\lambda=\infty$ and $EM_\infty\ge 0$.

The desired
statement holds true owing to $|M_{S_\lambda}-\lambda|\le
|\triangle M_{S_\lambda}|\le K$, that is,
$|M_{S_\lambda}-\lambda|,\ \lambda>0$ is a uniformly integrable
family. \qed

\subsection{Proof of part (iii)}

\subsubsection{\bf Auxiliary lemmas}

\begin{lemma}\label{lem-3.n1}
Under the assumptions of Theorem \ref{theo-n1} {\rm (}iii{\rm )},
$$
\lim_{\lambda\downarrow
0}E\frac{1}{\lambda}\Big(1-e^{-\log\mathscr{E}_\infty(\lambda
)}\Big) =EM_\infty.
$$
\end{lemma}
\begin{proof}
Recall that $\lambda\le\varepsilon$ for $\varepsilon$ involved in
assumption (ii). Since by Proposition \ref{pro-2.1} $\mathfrak{z}_t(\lambda)$
a uniformly integrable martingale, we have
$E\mathfrak{z}_\infty(\lambda )=1$. Hence,
$$
\begin{aligned}
E\frac{1}{\lambda}\Big(1- e^{-\log\mathscr{E}_\infty(\lambda
)}\Big)&= E\frac{1}{\lambda}\Big(\mathfrak{z}_\infty(\lambda )-
e^{-\log\mathscr{E}_\infty(\lambda )}\Big)
\\
&=E\frac{1}{\lambda}{\Big( e^{\lambda
M_\infty}-1\Big)}e^{-\log\mathscr{E}_\infty(\lambda )}.
\end{aligned}
$$
The required statement follows from the relation
$$
\begin{aligned} &
\lim_{\lambda\downarrow
0}\frac{1}{\lambda}e^{-\log\mathscr{E}_\infty(\lambda )}\Big(
e^{\lambda  M_\infty}-1\Big)=M_\infty,
\\
& \frac{1}{\lambda}e^{-\log\mathscr{E}_\infty(\lambda )}\big|
e^{\lambda  M_\infty}-1\big|\le e^{\varepsilon M_\infty}
\end{aligned}
$$
and the assumption $Ee^{\varepsilon M_\infty} <\infty $, see
\eqref{oxbla0}.
\end{proof}

\begin{lemma}\label{lem-3.n2}
Under the assumptions of Theorem \ref{theo-n1} {\rm (}iii{\rm )},

$$
\lim_{\lambda\downarrow
0}E\frac{1}{\lambda}\Big(1-e^{-\frac{\lambda^2}{2}\langle M
\rangle_\infty}\Big) =EM_\infty.
$$
\end{lemma}
\begin{proof}
Due to Lemma \ref{lem-3.n1}, suffice it to show that
\begin{equation}\label{3.n4}
\lim_{\lambda\downarrow
0}E\frac{1}{\lambda}\Big|e^{-\log\mathscr{E}_\infty(\lambda )}
-e^{-\frac{\lambda^2}{2}\langle M\rangle_\infty}\Big|=0.
\end{equation}
In order to verify \eqref{3.n4}, we estimate
$\log\mathscr{E}_\infty(\lambda)$ from above and below via
$\frac{\lambda^2}{2}\langle M\rangle_\infty$. Owing to $
\log\mathscr{E}_\infty(\lambda)\le \frac{\lambda^2}{2}\langle
M^c\rangle_\infty +G_\infty(\lambda), $ we have
\begin{equation}\label{new3.4}
\log\mathscr{E}_\infty(\lambda)\le\frac{\lambda^2}{2}\langle
M\rangle_\infty \Big[1+\frac{\lambda}{3}Ke^{\lambda K}\Big].
\end{equation}
Further, with
$$
G^c_\infty(\lambda)=\int_0^\infty\int_{|z|\le K}\big(e^{\lambda
z}-1-\lambda z\big) \nu^c(dt,dz),
$$
where $\nu^c(dt,dz):=\nu(dt,dz)-\nu(\{t\},dz)$, and
$
\Phi(\lambda,K)=1-\lambda Ke^{\lambda K},
$
we get
\begin{equation}\label{new3.5}
\begin{aligned}
\log\mathscr{E}_\infty(\lambda)
&=\frac{\lambda^2}{2}\langle M^c\rangle_\infty+
G^c_\infty(\lambda)+\sum_{t>0}\log\big(1+\triangle
G_t(\lambda)\big)
\\
&\ge\frac{\lambda^2}{2}\langle M^c\rangle_\infty+ \Phi(\lambda,K)
\int_0^\infty\int_{|z|\le K}\frac{\lambda^2}{2}
z^2\nu^c(dt,dz)
\\
&\quad
+\sum_{t>0}\log\Bigg(1+\Phi(\lambda,K)
\int_{|z|\le K}\frac{\lambda^2}{2}z^2\nu(\{t\},dz)\Bigg).
\end{aligned}
\end{equation}
We choose $\lambda$ so small to have $1-\lambda Ke^{\lambda
K}>0$ and estimate from below the ``$\sum_{t>0}\log$'' in the last
line from the above inequality by applying
$$
\log(1+x)\ge x-\frac{1}{2}x^2, \ x\ge 0.
$$
This gives us the bound
\begin{gather*}
\sum_{t>0}\log\Bigg(1+\Phi(\lambda,K)
\int_{|z|\le K}\frac{\lambda^2}{2}z^2\nu(\{t\},dz)\Bigg)
\\
\ge \Phi(\lambda,K) \int_{|z|\le
K}\frac{\lambda^2}{2}z^2\nu(\{t\},dz)
-\frac{1}{2}\Phi^2(\lambda,K) \left(\int_{|z|\le
K}\frac{\lambda^2}{2}z^2\nu(\{t\},dz)\right)^2.
\end{gather*}
Since $\nu(\{t\},|z|\le K)\le 1$, by the Cauchy--Schwarz
inequality we find that
\begin{gather*}
\left(\int_{|z|\le K}\frac{\lambda^2}{2}z^2\nu(\{t\},dz)\right)^2
\\
\le \frac{\lambda^4}{4}\int_{|z|\le K}z^4\nu(\{t\},dz) \le
\frac{\lambda^4K^2}{4}\int_{|z|\le K}z^2\nu(\{t\},dz).
\end{gather*}
So, finally we get
\begin{gather}\label{new3.6}
\sum_{t>0}\log\Bigg(1+\Phi(\lambda,K)
\int_{|z|\le K}\frac{\lambda^2}{2}z^2\nu(\{t\},dz)\Bigg)
\nonumber\\
\quad\ge \left(\Phi(\lambda,K)
-\frac{\lambda^2}{8}K^2\Phi^2(\lambda,K)\right)
\int_{|z|\le K}\frac{\lambda^2}{2}z^2\nu(\{t\},dz)
\end{gather}
and now choose $\lambda$ so small to have
\begin{equation}\label{new3.7}
\Phi(\lambda,K)
-\frac{\lambda^2}{8}K^2\Phi^2(\lambda,K)\ge
1-\lambda C>0
\end{equation}
for some constant $C>0$. Combining now \eqref{new3.5},
\eqref{new3.6} and \eqref{new3.7}, we may choose a generic
positive constant $C$ and sufficiently small $\lambda$ such that
$\mathscr{E}_\infty(\lambda)\ge
\big[1-C\lambda\big]\frac{\lambda^2}{2}\langle M\rangle_\infty. $
Hence and with \eqref{new3.4}, for some generic positive constant
$C>0$ and sufficiently small $\lambda>0$ we have
$$
0<\big[1-C\lambda\big]\frac{\lambda^2}{2}\langle M\rangle_\infty
\le\log\mathscr{E}_\infty(\lambda)\le
\big[1+C\lambda\big]\frac{\lambda^2}{2} \langle M\rangle_\infty.
$$
These inequalities provide
$$
\frac{1}{\lambda}\Big|e^{-\log\mathscr{E}_\infty(\lambda )}
-e^{-\frac{\lambda^2}{2}\langle M\rangle_\infty}\Big|\le
C\frac{\lambda^2}{2}\langle M\rangle_\infty
e^{-\frac{\lambda^2}{2}\langle M\rangle_\infty}\xrightarrow[\lambda\to 0]{}0.
$$
Since $xe^{-x}\le e^{-1}$, the desired result holds by Lebesgue's
dominated theorem.
\end{proof}

\begin{lemma}\label{lem-3.n3}
Under the assumptions of Theorem \ref{theo-n1} {\rm (}iii{\rm )},
$$
\lim\limits_{\lambda\to\infty}\lambda P\big(\langle
M\rangle^{1/2}_\infty>\lambda\big)=c \Leftrightarrow
\lim_{\lambda\to\infty}\lambda
P\big([M,M]^{1/2}_\infty>\lambda\big)=c.
$$
\end{lemma}
\begin{proof}
It suffices to establish
\begin{equation}\label{4.1aa}
\begin{split}
& \varlimsup_{\lambda\to
0}\frac{P\big([M,M]^{1/2}_\infty>\lambda\big)} {P\big(\langle
M\rangle^{1/2}_\infty>\lambda\big)}\le 1,
\\
& \varliminf_{\lambda\to
0}\frac{P\big([M,M]^{1/2}_\infty>\lambda\big)} {P\big(\langle
M\rangle^{1/2}_\infty>\lambda\big)}\ge 1.
\end{split}
\end{equation}
Set $L=[M,M]-\langle M\rangle$. Since $ [M,M]_\infty\le \langle
M\rangle_\infty+\sup_{t\ge 0}|L_t|$,  applying the elementary
inequality $ (c+d)^{1/2}\le c^{1/2}+d^{1/2}$, we find that
\begin{align}\label{3.n7}
& P\big([M,M]^{1/2}_\infty>\lambda\big)\le P\big([\langle
M\rangle_\infty+\sup_{t\ge 0} |L_t|]^{1/2}>\lambda\big)
\nonumber\\
&\quad \le P\big(\langle M\rangle^{1/2}_\infty+\sup_{t\ge 0}
|L_t|^{1/2}>\lambda\big)
\nonumber\\
&\quad \le P\big(\langle M\rangle^{1/2}_\infty>(1-a)\lambda\big)+
P\big(\sup_{t\ge 0}|L_t|>a\lambda\big), \ a\in(0,1).
\end{align}
With $\lambda_a=(1-a)\lambda$, the resulting bound can be
rewritten as:
\begin{equation}\label{eb}
\lambda P\big([M,M]^{1/2}_\infty>\lambda\big)\le
(1-a)^{-1}\lambda_aP\big(\langle
M\rangle^{1/2}_\infty>\lambda_a\big) +\lambda P\big(\sup_{t\ge
0}|L_t|^{1/2}>a\lambda\big).
\end{equation}
So, we shall deal with the evaluation from above of
$P\big(\sup_{t\ge 0}|L_t|^{1/2}>a\lambda\big)$.
A helpful tool here is the inequality: for some  absolute positive constant $C$,
any stopping time $\tau$ and $K$ being a bound for $|\triangle M|$,
\begin{equation}\label{3.9z}
E\sup_{t\le\tau}|L_t|^2\le CK^2 E\langle M\rangle_\tau.
\end{equation}
In order to establish \eqref{3.9z}, we use the following facts:

- $L$ is the purely discontinuous local martingale with
\begin{align*}
[L,L]_t&=\sum_{s\le t}(\triangle L_s)^2 =\sum_{s\le
t}\big((\triangle M_s)^2-\triangle\langle M\rangle_s\big)^2
\\
&=\sum_{s\le t}\left(\int_{|z|\le
K}z^2(\mu(\{s\},dz)-\nu(\{s\},dz)\right)^2,
\end{align*}

-
$
\langle L\rangle_t=\int_0^t\int_{|z|\le
K}z^4(\nu(ds,dz)- \sum_{s\le t}\Big(\int_{|z|\le
K}z^2\nu(\{s\},dz)\Big)^2,
$

-
$
\begin{aligned}
\langle L\rangle_t\le \int_0^t\int_{|z|\le
K}z^4\nu(ds,dz) \le K^2\int_0^t\int_{|z|\le
K}z^2\nu(\{ds,dz)\le K^2\langle M\rangle_t,
\end{aligned}
$

- $K^2\langle M\rangle-\langle L\rangle$ is the increasing process.

\medskip
\noindent
Now, we refer to the Burkholder-Gundy inequality (see e.g. Theorem
1.9.7 in \cite{LSMar}): for any stopping time $\tau$,
$$
E\sup_{t\le\tau}|L_t|^2\le CE[L,L]_\tau.
$$

Due to the relations
$E[L,L]_\tau=E\langle L\rangle_\tau$ and $K^2\langle
M\rangle_\tau\ge \langle L\rangle_\tau$ (recall that $K^2\langle
M\rangle\ge \langle L\rangle$), we have
$
E\langle L\rangle_\tau\le K^2 E\langle M\rangle_\tau,
$
that is, \eqref{3.9z} is valid.
Due to \eqref{3.9z} and the fact that $\langle M\rangle$ is a
predictable process, the Lenglart--Rebolledo inequality (see,
e.g., Theorem 1.9.3 in \cite{LSMar}) is  applicable (notice that
$\{\sup_{t\ge 0}|L_t|^{1/2}>a\lambda\}\equiv\{\sup_{t\ge
0}|L_t|>a^2\lambda^2\}$), so that,
\begin{gather}\label{3.n8}
P\Big(\sup_{t\ge 0}|L_t|^{1/2}>a\lambda\big)\le
\frac{\lambda^{5/2}}{a^4\lambda^4}+P\big(CK^2\langle
M\rangle_\infty> \lambda^{5/2}\Big)
\nonumber\\
=\frac{\lambda^{5/2}}{a^4\lambda^4}+P\Big(\langle
M\rangle^{1/2}_\infty> \lambda^{5/4}/(C^{1/2}K)\Big).
\end{gather}
Hence, with $r=1/(C^{1/2}K)$ and $\lambda_r=r\lambda^{5/4}$,
\begin{equation}\label{4.4a}
\lambda P\Big(\sup_{t\le\mbox{}_{T_x}}|L_t|^{1/2}>a\lambda\Big)
\le \frac{1}{a^4\lambda^{1/2}}+\frac{1}{r\lambda^{1/4}}\lambda_r
P\Big(\langle M\rangle^{1/2}_\infty> \lambda_r\Big).
\end{equation}
Now, \eqref{eb} and \eqref{4.4a} provide
\begin{gather*}
\lambda P\Big([M,M]^{1/2}_\infty>\lambda\big)
\\
\le(1-a)^{-1}\lambda_a P\big(\langle
M\rangle^{1/2}_\infty>\lambda_a\Big)
+
\frac{1}{a^4\lambda^{1/2}}+\frac{r}{\lambda^{1/4}}\lambda_r
P\big(\langle M\rangle^{1/2}_\infty>\lambda_r\big).
\end{gather*}

\medskip
Assume that $c>0$. Then, we get
\begin{gather*}
\frac{P\big([M,M]^{1/2}_\infty>\lambda\big)} {P\big(\langle
M\rangle^{1/2}_\infty>\lambda\big)} \le
\frac{(1-a)^{-1}\lambda_aP\big(\langle
M\rangle^{1/2}_\infty>\lambda_a\big) }{ \lambda P\big(\langle
M\rangle^{1/2}_\infty>\lambda\big)}
\\
+\frac{\frac{1}{a^4\lambda^{1/2}}
+\frac{r}{\lambda^{1/4}}\lambda_rP\big(\langle
M\rangle^{1/2}_\infty> \lambda_r\big)}{ \lambda P\big(\langle
M\rangle^{1/2}_\infty>\lambda\big)}\xrightarrow[\lambda\to \infty]{}
\frac{1}{1-a}\xrightarrow[a\to 0]{}1
\end{gather*}
and the first part from \eqref{4.1aa}.

Since the second part from
\eqref{4.1aa} is established similarly, we give only a sketch of
the proof. The use of
\[
P\big(\langle M\rangle^{1/2}>\lambda\big)\le
P\big([M,M]^{1/2}>(1-a)\lambda\big)+ P\big(\sup_{t\ge
0}|L_t|>a\lambda\big), \ a\in(0,1)
\]
provides
$$
\frac{P\big([M,M]\rangle^{1/2}_\infty>(1-a)\lambda\big)}
{P\big(\langle M\rangle^{1/2}_\infty>\lambda\big)}\ge 1-
\frac{P\big(\sup_{t\ge 0}|L_t|>a\lambda\big)} {P\big(\langle
M\rangle^{1/2}_\infty>\lambda\big)}
$$
and the result.

\medskip
If $c=0$, we replace $M$ by $M+M'$, where $M'$ is independent of
$M^c$ local continuous  martingale with $M'_0=0$ and $\langle
M'\rangle_\infty<\infty$ a.s. and
$$
\lim_{\lambda\to\infty}\lambda P\big(\langle
M'\rangle^{1/2}_\infty>\lambda\big)=c'>0.
$$
Now, taking into account the  obvious relations
\begin{align*}
& [M+M',M+M']=[M,M]+[M',M']\ {\rm and} \ \langle
M+M'\rangle=\langle M\rangle+\langle M'\rangle,
\end{align*}
with  $\delta\ne 0$ we find that $
\lim\limits_{\lambda\to\infty}\lambda P\big(\langle M+\delta
M'\rangle^{1/2}_\infty> \lambda\big)=\delta^2c'>0. $
So, by using the result already proved, we have
$$
\lim_{\lambda\to\infty}\lambda P\big([M+\delta M',M+\delta
M']^{1/2}_\infty
>\lambda\big)=\delta c'
$$
and so, by $ P\big([M+\delta M',M+\delta
M']^{1/2}_\infty>\lambda\big) \ge
P\big([M,M]^{1/2}_\infty>\lambda\big)$,
we find that
$$
\varlimsup\limits_{\lambda\to 0}\lambda
P\big([M,M]^{1/2}_\infty>\lambda\big) \le \delta
c'\xrightarrow[\delta\to 0]{}0.
$$
\end{proof}

\subsubsection{\bf Final part of the proof for (iii)}

We refer to the Tauberian theorem.

\medskip
\noindent {\bf Theorem**.} \label{Tauber}{\rm (Feller,
\cite{Feller}, XIII.5, Example (c))} {\it Let $X$ be a nonnegative
random variable such that $\lim\limits_{\lambda\downarrow
0}\frac{1}{\lambda}\Big(1-Ee^{-\frac{\lambda^2}{2} X}\Big)$ exists
in $\mathbb{R}$, then
$$
\sqrt{\frac{2}{\pi}}\lim_{\lambda\downarrow 0}\frac{1}{\lambda}
\Big(1-Ee^{-\frac{\lambda^2}{2} X}\Big)=
\lim_{\lambda\to\infty}\lambda P(X^{1/2}>\lambda).
$$
}

\bigskip
Now, we are in the position to finish the proof of (ii). Letting
$X=\langle M\rangle_\infty$, we find that
$$
\sqrt{\frac{2}{\pi}}\lim_{\lambda\downarrow
0}\frac{1}{\lambda}\Big(1-Ee^{-\frac{\lambda^2}{2}\langle
M\rangle_\infty}\Big)= \lim_{\lambda\to\infty}\lambda P(\langle
M\rangle^{1/2}_\infty>\lambda).
$$
At the same time, Lemmas \ref{lem-3.n1} and \ref{lem-3.n2} provide
$$
\lim_{\lambda\downarrow
0}\frac{1}{\lambda}\Big(1-Ee^{-\frac{\lambda^2}{2}\langle
M\rangle_\infty}\Big)= \sqrt{\frac{2}{\pi}}EM_\infty
$$
while by Lemma \ref{lem-3.n3}
$\lim\limits_{\lambda\to\infty}\lambda P\big([M,M]^{1/2}_\infty>
\lambda\big)=\sqrt{\frac{2}{\pi}}EM_\infty$. \qed

\subsection{Supplement}

As it was mentioned in Introduction,  the condition $|\triangle
M|\le K$ might be too restrictive to be valid for serving some examples. It is known
from \cite{Nov} that this condition can be replaced by a weaker
one and so more useful for applications. An analog of this result is given
below.

\begin{theorem}\label{theo-n3}
Let $M\in\mathscr{M}^2_{\rm loc}, \langle M\rangle_\infty<\infty$
a.s., $M^+ \in \mathcal{D}$ and \eqref{oxbla0} holds. Assume also
that  there exist nonnegative integrable random variables
$\zeta_1$, $\zeta_2$ such that for all sufficiently small
$\lambda>0$

\begin{equation}\label{rest}
\frac{\lambda^2}{2}\langle
M\rangle_\infty(1-|\lambda|\zeta_1)^+\le
\log\mathscr{E}_\infty(\lambda)\le \frac{\lambda^2}{2}\langle
M\rangle_\infty(1+|\lambda|\zeta_2).
\end{equation}
Then
$$
\lim\limits_{\lambda\to\infty}\lambda P\big(\langle
M\rangle^{1/2}_\infty > \lambda\big)
=\sqrt{\frac{2}{\pi}}EM_\infty.
$$
\end{theorem}

\begin{proof}
Notice that only \eqref{3.n4} has to be verified under
\eqref{rest}.

By \eqref{rest}, we have
$$
\begin{aligned}
\frac{1}{\lambda}\Big|e^{-\log\mathscr{E}_\infty(\lambda )}
-e^{-\frac{\lambda^2}{2}\langle M\rangle_\infty}\Big|&\le
\Big(\zeta_2\vee\frac{|1-(1-\zeta_1\lambda)^+|}{\lambda}\Big)\frac{\lambda^2}{2}\langle
M\rangle_\infty e^{-\frac{\lambda^2}{2}\langle M\rangle_\infty}
\\
& \le \big(\zeta_2\vee\zeta_1\big)\frac{\lambda^2}{2}\langle
M\rangle_\infty e^{-\frac{\lambda^2}{2}\langle M\rangle_\infty}.
\end{aligned}
$$
The right-hand side of this inequality converges to zero, as
$\lambda\to 0$, and is bounded by $e^{-1}(\zeta_2\vee\zeta_1)$.
Hence, in order to get \eqref{3.n4} suffices it  to allude on the Lebesgue
dominated convergence theorem.
\end{proof}

\medskip
{\bf Acknowledgements.} The authors gratefully acknowledge their
colleagues J. Stoyanov,  E. Shinjikashvili and anonymous reviewers for
comments improving presentation of the material.


\end{document}